\documentclass[a4paper, 11pt]{article}
\usepackage{graphicx,tcolorbox}
\usepackage{booktabs,caption,longtable,multicol,multirow}
\usepackage{subcaption,lscape}
\usepackage[linesnumbered,boxruled,lined,commentsnumbered]{algorithm2e}
\usepackage[colorlinks,citecolor=purple,linkcolor=blue]{hyperref}
\usepackage[margin=1.2in]{geometry}
\usepackage[numbers]{natbib}
\newcommand{\cupdlp}{\textrm{cuPDLP-C}}
\newcommand{\cupdlpjl}{\textrm{cuPDLP.jl}}
\newcommand{\copt}{\textrm{COPT}}

\usepackage{amsmath,amssymb,mathrsfs,amsthm}
\DeclareMathOperator*{\proj}{\textbf{proj}}

\title{cuPDLP-C: A Strengthened Implementation of cuPDLP for Linear Programming by C language}

\author{
Haihao Lu\textsuperscript{1}\thanks{Contact: \texttt{haihao.lu@chicagobooth.edu}} \and Jinwen Yang\textsuperscript{2} \and Haodong Hu\textsuperscript{3} \and Qi Huangfu\textsuperscript{4} \and Jinsong Liu\textsuperscript{3} \and Tianhao Liu\textsuperscript{3} \and Yinyu Ye\textsuperscript{5} \and Chuwen Zhang\textsuperscript{3} \and Dongdong Ge\textsuperscript{3}\thanks{Correspondence to: \texttt{ge.dongdong@mail.shufe.edu.cn}. This research is partially supported by the National Natural Science Foundation of China (NSFC) [Grant NSFC-72150001, 72225009, 72394360, 72394365].}\\\\
\textsuperscript{1}The University of Chicago, Booth School of Business \\
\textsuperscript{2}The University of Chicago, Department of Statistics \\
\textsuperscript{3}Shanghai University of Finance and Economics \\
\textsuperscript{4}Cardinal Operations \quad
\textsuperscript{5}Stanford University
}
\date{\today}

\begin{document}
\maketitle
\begin{abstract}
    A recent GPU implementation of the Restarted Primal-Dual Hybrid Gradient Method for Linear Programming (\cupdlpjl{}) was proposed in \cite{lu_cupdlpjl_2023}. Its computational results demonstrate the significant computational advantages of the GPU-based first-order algorithm on certain large-scale problems. The average performance also achieves a level close to commercial solvers for the first time in history. However, due to limitations in experimental hardware and the disadvantage of implementing the algorithm in Julia compared to C language, neither the commercial solver nor cuPDLP reached their maximum efficiency. Therefore, in this report, we have re-implemented and optimized cuPDLP in C language. Utilizing state-of-the-art CPU and GPU hardware, we extensively compare cuPDLP with the best commercial solvers. The experiments further highlight its substantial computational advantages and potential for solving large-scale linear programming problems. We also discuss the profound impact this breakthrough may have on mathematical programming research and the entire operations research community.
\end{abstract}

\section{Introduction}
Recently, there has been heated interest in using first-order methods for solving linear programs that use matrix-vector multiplication instead of factorizations as the major workhorse. Due to the success of primal-dual hybrid gradient on linear programming (PDLP) \cite{applegate_practical_2022}, exciting progress has been made from both the theoretical and computational aspects. As the original analysis of PDLP \cite{applegate_faster_2023} involves a condition number that is difficult for evaluation, a vein of research is devoted to improving the convergence analysis from a geometric perspective and special type problems; see \cite{lu_geometry_2023,hinder_worst-case_2023,xiong_computational_2024} for example. 

On the practical side, the initial Julia (CPU) version has been reimplemented in C++ as part of the linear programming kit in Google OR Tools. \citet{lu_cupdlpjl_2023} introduced a GPU-based PDLP that provided an affirmative answer to a long-standing question of \emph{whether a Linear Programming solver can capitalize on GPU-type architectures to enable efficient parallel computing}. Although a rigorous implementation was provided in \cite{lu_cupdlpjl_2023} in Julia programming language, due to the limitation of computing resources, a few questions persist on whether (a) an optimized C/C++ implementation and (b) more advanced hardware can further unleash the potential of the proposed method.

The goal of this paper is to complement the results in \citet{lu_cupdlpjl_2023} by providing a strengthened implementation of cuPDLP in C programming language with generic CUDA C/C++ API under the name \cupdlp{}. Furthermore, to fully demonstrate the pros and cons of the method, we compare \cupdlp{} with the state-of-the-art interior-point-based optimizer COPT \cite{ge_cardinal_2022} under the best available computing hardware. We run COPT on AMD Ryzen 9 5900X with 512GB RAM, whereas the GPU-based \cupdlp{} is tested on NVIDIA H100 80GB. Our implementation is public at \texttt{https://github.com/COPT-Public/cuPDLP-C}.

\section{Overview of the Method}
We give a minimal review of the method and the techniques used in our implementation.

\paragraph{Linear programming.} We consider the following primal-dual linear programming (LP) problems
\begin{equation}
    \begin{aligned}
        \min\ & c^\top x \\
        \text{s.t.}\ & Ax = b \\
        & Gx \geq h \\
        & l \leq x \leq u
    \end{aligned}
    \quad\quad\quad\quad
    \begin{aligned}
        \max\ & q^\top y + l^\top \lambda^+ - u^\top \lambda^- \\
        \text{s.t.}\ & c - K^\top y = \lambda \\
        & y_{m_1+1:m_1+m_2} \geq 0 \\
        & \lambda \in \Lambda
    \end{aligned}
    \label{eq:LP}
\end{equation}
where $A\in\mathbb{R}^{m_1\times n}$, $G\in\mathbb{R}^{m_2\times n}$, $c\in\mathbb{R}^{n}$, $b\in\mathbb{R}^{m_1}$, $h\in\mathbb{R}^{m_2}$, $l\in(\mathbb{R}\cup\{-\infty\})^{n}$, $u\in(\mathbb{R}\cup\{+\infty\})^{n}$, $K^\top = (A^\top, G^\top)$, $q^\top = (b^\top, h^\top)$, and
\begin{equation*}
    \Lambda = \Lambda_1\times \cdots \times \Lambda_n,\ \Lambda_i = \begin{cases}
        \{0\} & l_i = - \infty, u_i = +\infty \\
        \mathbb{R}^- & l_i = - \infty, u_i \in \mathbb{R} \\
        \mathbb{R}^+ & l_i \in \mathbb{R}, l_i = + \infty \\
        \mathbb{R} & \text{otherwise}
    \end{cases}
\end{equation*}

\paragraph{Optimality termination criteria.}
cuPDLP terminates when the primal-dual solution $(x,y,\lambda)$ satisfies bound constraints and
\begin{itemize}
    \item Primal feasibility
    \begin{equation*}
        \left\| \begin{matrix}
            Ax - b \\ [h - Gx]^+
        \end{matrix} \right\|_2 \leq \epsilon (1 + \|q\|_2)
    \end{equation*}
    \item Dual feasibility
    \begin{equation*}
        \left\| c - K^\top y - \lambda \right\|_2 \leq \epsilon (1 + \|c\|_2)
    \end{equation*}
    \item Duality gap
    \begin{equation*}
        \left|q^\top y + l^\top \lambda^+ - u^\top \lambda^- - c^\top x\right| \leq \epsilon (1 + |q^\top y + l^\top \lambda^+ - u^\top \lambda^-| + |c^\top x|)
    \end{equation*}
\end{itemize}
Notice that dual variable $\lambda$ for the bound constraints is not explicitly generated during cuPDLP; we compute $\lambda = \proj_{\Lambda}(c - K^\top y)$.

\paragraph{Restarted PDHG.} To solve \eqref{eq:LP}, we can equivalently solve the saddle-point problem
\begin{equation}\label{eq:saddle-point}
    \min_{x\in X}\max_{y\in Y} \ \mathcal{L}(x,y) = c^\top x - y^\top Kx + q^\top y 
\end{equation}
with $X=\{x\in\mathbb{R}^n: l\leq x \leq u\}$ and $Y=\mathbb{R}^{m_1}\times\mathbb{R}^{m_2}_+$.

Using restarted PDHG to solve \eqref{eq:saddle-point}, the $n$th inner loop of the algorithm starts from $(x^{n,0}, y^{n,0})$ and consists of three parts
\begin{enumerate}
    \item Primal update
    \begin{equation*}
        x^{n,t+1} = \proj_{X}(x^{n,t} - \frac{\eta}{\omega^n} (c - K^\top y^{n,t}))
    \end{equation*}
    \item Dual update
    \begin{equation*}
        y^{n,t+1} = \proj_{Y}(y^{n,t} + \eta\omega^n (q - K(2x^{n,t+1} - x^{n,t})))
    \end{equation*}
    \item Restart
    \item[] We keep track of the averaged primal-dual solution, and restart the $(n+1)$th primal-dual update loop from the averaged solution if certain restart conditions hold.
\end{enumerate}

\paragraph{Restart criteria.}
For the GPU version, restart criteria are slightly different from the one in \cite{applegate_faster_2021, applegate_practical_2022}. In \cite{applegate_faster_2021}, it is proved that restart can improve the convergence rate of PDHG to optimum for solving LP. However, in \cite{applegate_faster_2021}, checking restart criteria requires computing a trust region subproblem, which is not friendly for GPU. Therefore, cuPDLP maintains the restart framework but replaces the subproblem with the GPU-friendly KKT error
\begin{equation}
    \text{KKT}_\omega(z) = \sqrt{\omega^2\left\|\begin{matrix}
        Ax - b \\ [h - Gx]^+
    \end{matrix}\right\|_2^2 + \frac{1}{\omega^2}\left\| c - K^\top y - \lambda \right\|_2^2 + (q^\top y + l^\top\lambda^+ - u^\top \lambda^- - c^\top x)^2}
\end{equation}
where $z = (x,y)$.

To be more specific, after choosing the restart candidate
\begin{equation*}
    z^{n, t+1}_c = \begin{cases}
        z^{n, t+1} & \text{KKT}_{\omega^n}(z^{n,t+1}) < \text{KKT}_{\omega^n}(\bar{z}^{n,t+1}) \\
        \bar{z}^{n,t+1} & \text{otherwise}
    \end{cases}
\end{equation*}
cuPDLP will restart if one of the three conditions holds:
\begin{enumerate}
    \item Sufficient decay in KKT error
    \begin{equation*}
        \text{KKT}_{\omega^n}(z^{n,t+1}_c) \leq 0.2 \text{KKT}_{\omega^n}(z^{n,0})
    \end{equation*}
    \item Necessary decay and no local progress in KKT error
    \begin{equation*}
        \text{KKT}_{\omega^n}(z^{n,t+1}_c) \leq 0.8 \text{KKT}_{\omega^n}(z^{n,0})\ \text{and}\ \text{KKT}_{\omega^n}(z^{n,t+1}_c) > \text{KKT}_{\omega^n}(z^{n,t}_c)
    \end{equation*}
    \item Long inner loop
    \begin{equation*}
        t \geq 0.36 k
    \end{equation*}
    where $k$ is the total iteration number.
\end{enumerate}

\section{Numerical Experiments}
In this section, we provide numerical results of \cupdlp{} and the comparisons to COPT. 

\subsection{Experimental Setup}

\paragraph{Hardware and software.}
The original experiments in \cite{lu_cupdlpjl_2023} were conducted in NVIDIA Tesla V100-PCIe-16GB and Intel Xeon Gold 6248R CPU 3.00GHz for GPU-based and CPU-based methods, respectively. To enable comprehensive understanding for \cupdlp{} under advanced hardware, we run the CPU solvers on AMD Ryzen 9 5900X, whereas the GPU ones are tested on NVIDIA H100 80GB HBM3. 
Except for techniques like scaling, we also use different presolving modules in cuPDLP-C. 
We directly reproduce the experiments of \cupdlpjl{} with the hardware at hand; the results (see \autoref{fig.benchmark.time}) show the C implementation has a clear advantage over the Julia version. For such reasons, we focus on the study of \cupdlp{} without further investigation of \cupdlpjl{}.

\paragraph{Benchmark datasets.}

We intensively conduct the performance experiments on both classical benchmark sets and huge-size instances. Firstly, we test LP solvers on two classic benchmark datasets, including MIPLIB 2017 subset and Mittelmann’s LP benchmark set. The MIPLIB 2017 collection is a classic mixed-integer programming collection. We select 383 instances based on the same criteria in \cite{lu_cupdlpjl_2023} and solve their LP relaxation. Mittelmann's LP benchmark is a famous benchmark containing several large-scale challenging LP instances, among which we utilize 49 public instances.  Secondly, to further probe the capability boundary of cuPDLP, we collect several huge-sized LP problems from different industries, such as PageRank, 
% unit commitment, 
supply chain management and so on. As interior-point and simplex methods heavily rely on matrix factorization, which is computationally infeasible at this magnitude, cuPDLP currently stands as the only realistic method for finding nearly optimal solutions in such challenging scenarios.

The running time per instance is delayed in \autoref{fig.mittelman.time} and \autoref{fig.mip383.time}. Results on the huge-size LP collection are reported in \autoref{fig.selected.time}. For Mittelmann's dataset and the MIPLIB collection, we present the overall performance on two benchmarks in \autoref{fig.benchmark.time}.

% \paragraph{A collection of huge-size LP problems}
% To demonstrate the effectiveness of both methods on huge-size problems, we collect a selective list of huge-size linear programming problems. We report the running time in \autoref{fig.selected.time}.

\paragraph{Evaluation.} We evaluate the solvers by shifted geometric mean (SGM) of solving time. For $n$ instances, SGM is defined as $(\prod_{i=1}^{n}(t_i + \Delta))^{1/n} - \Delta$ where $t_i$ is the solving time for the $i$th instance. If the instance is unsolved, its solving time is set to the time limit. We set $\Delta=10$ and denote the metric as SGM10. For fair comparisons, we exclude the running time needed for different presolving and scaling techniques. 

\subsection{MIPLIB relaxation and Mittelmann's LP benchmark}
We test \cupdlp{} with various presolvers, \cupdlpjl{} under the tolerance of $10^{-4}$ and $10^{-8}$, and \copt{} under the tolerance of $10^{-8}$.  As shown in \autoref{fig.mip383.time}, COPT can solve all 383 instances, while \cupdlp{} solves 379 under $10^{-4}$ tolerance.
% \subsection{}
% We test \cupdlp{} with various presolvers under the tolerance of $10^{-4}$ and $10^{-8}$, and \copt{} under the tolerance of $10^{-8}$. 
It is observed in \autoref{fig.mittelman.time} that again COPT solves the most instances and thus achieves the smallest SGM10. \cupdlp{} performs around 2 to 4 times slower than COPT with different presolvers under $10^{-4}$ accuracy.

For internal comparison, \cupdlp{} significantly outperforms \cupdlpjl{}, which reveals that reimplementing cuPDLP from Julia in C brings about a 50\% increase in speed. Furthermore, our experiments also indicate presolvers can further speed up the convergence. With the presolving modules in CLP and COPT, \cupdlp{} demonstrates itself as a highly competitive LP solver, which is readily faster than most of the open-source LP solvers according to \cite{mittelmann_lp_2023}.

\subsection{Some Large Problems}

% \begin{table}
%     \small
%     \centering
%     \caption{Overall performance (SGM10) on different datasets}\label{fig.benchmark.time}
%     \begin{tabular}{l|r|r|r|r|r|r|r}
%         \toprule
%         \multirow{3}{*}{Dataset} & \copt      & \multicolumn{6}{|c}{\cupdlp{}}                                                                                                   \\
%                                  &            & \multicolumn{2}{|c}{No Presolve} & \multicolumn{2}{|c}{COPT} & \multicolumn{2}{|c}{CLP}                                          \\
%                                  & $10^{-8}$  & $10^{-4}$                        & $10^{-8}$                 & $10^{-4}$                & $10^{-8}$   & $10^{-4}$  & $10^{-8}$   \\
%         \midrule
%         Mittelmann (49)          & 13.81 (48) & 57.54 (43)                       & 172.98 (39)               & 57.54 (40)               & 173.17 (35) & 40.87 (43) & 143.13 (37) \\
%         MIPLIB (383)             & 3.11 (383) & 28.37                            & 10.26                                                                                         \\
%     \end{tabular}
% \end{table}
\begin{table}
    \small
    \centering
    \caption{Overall performance (SGM10) on different datasets}\label{fig.benchmark.time}
    \begin{tabular}{l|r|r|r|r|r}
        \toprule
        Dataset                          & Optimizer                  & Presolver             & Tol.      & SGM10 & Solved \\
        \midrule
        \multirow{7}{*}{MIPLIB (383) }   & COPT                       &               -        & $10^{-8}$ & 3.11 &383         \\
        \cmidrule{2-6}
                                         & \multirow{6}{*}{\cupdlp{}} & \multirow{2}{*}{COPT} & $10^{-4}$ & 5.43 &379         \\
                                         &                            &                       & $10^{-8}$ & 18.53 &369        \\
        \cmidrule{3-6}
                                         &                            & \multirow{2}{*}{CLP}  & $10^{-4}$ & 7.95 &372          \\
                                         &                            &                       & $10^{-8}$ & 21.89 & 362         \\
        \cmidrule{3-6}
                                         &                            & \multirow{2}{*}{No Presolve}    & $10^{-4}$ & 10.28 &370	          \\
                                         &                            &                       & $10^{-8}$ & 27.15 &359         \\
        \cmidrule{2-6}
                                         &  \multirow{2}{*}{\cupdlpjl{}} & \multirow{2}{*}{No Presolve} & $10^{-4}$ & 17.49 &370 \\
                                         &        &   & $10^{-8}$ & 35.69  & 355                   \\
        \midrule
        \multirow{7}{*}{Mittelmann (49)} & COPT                       &           -            & $10^{-8}$ & 13.81 &48          \\
        \cmidrule{2-6}
                                         & \multirow{6}{*}{\cupdlp{}} & \multirow{2}{*}{COPT} & $10^{-4}$ & 25.29 &46        \\
                                         &                            &                       & $10^{-8}$ & 110.22 &41         \\
        \cmidrule{3-6}
                                         &                            & \multirow{2}{*}{CLP}  & $10^{-4}$ & 33.97 & 45          \\
                                         &                            &                       & $10^{-8}$ & 125.95 & 38         \\
        \cmidrule{3-6}
                                         &                            & \multirow{2}{*}{No Presolve}    & $10^{-4}$ & 57.54 &43          \\
                                         &                            &                       & $10^{-8}$ & 172.98 &39         \\
        \bottomrule
    \end{tabular}
    \normalsize
\end{table}

% \begin{table}
%     \small
%     \centering
%     \caption{Performance on selected problems}\label{fig.selected.time}
%     \begin{tabular}{l|r|r}
%         \toprule
%         Instance       & Optimizer                  & Time. \\
%         \midrules
%       \bottomrule
%     \end{tabular}
% \end{table}   

To highlight the advantages of \cupdlp{} over COPT on extremely large instances, we pick a handful of representative instances, and then use different methods to find approximate primal-dual solutions. The target accuracy is set to $10^{-6}$ by default. The results are presented in \autoref{fig.selected.time}.

We include the instance ``zib03'' mentioned in \citet{kochProgressMathematicalProgramming2022}. The instance was made public in 2008 and first solved in April 2009, taking almost 139 days by CPLEX. A basis solution was only recently discovered in 2021.  

Another problem class arises from the LP formulation of the PageRank problem \citep{nesterov2014subgradient}. We generated PageRank instances from both randomly generated graphs \citep{applegate_practical_2022} and a real-world graph dataset \citep{snapnets}. In particular, we created two instances with $1$ million and $10$ million nodes, respectively. Additionally, we selected several instances with millions of nodes, resulting in extremely large LP instances with millions of rows and columns and tens of millions of non-zero coefficients, which are difficult for traditional LP solvers. 

We further introduce a set of supply-chain optimization problems from industrial applications at Cardinal Operations. The instances are modeled from an inventory management problem with transshipment decisions on a large-scale supply-chain network in a finite horizon $T$. We give several instances that stand for $T=10, 20, 40, 60$, respectively. The optimization problem produces plans for 600 facilities and customers with 200 commodities. For these problems, we set the gap tolerance to $10^{-5}$.

We generate instances of quadratic assignment problems (QAPs) using Adams-Johnson linearization \cite{adams_improved_1994} from QAPLIB \cite{burkard1997qaplib}. Each original QAP instance contains 50 locations. For the LP relaxation of Adams and Johnson linearization, we set the tolerance to $10^{-6}$ and the time limit to 3600s. Both simplex and barrier methods of COPT fail to provide solutions, while \cupdlp{} efficiently solves all five relaxations.

\begin{table}[htbp]
\small
  \centering
    \caption{Performance on selected problems}\label{fig.selected.time}
    \begin{tabular}{c|c|rrrrr}
    \toprule
    Source  & Instance & $m$     & $n$     & $nnz$   & COPT & \cupdlp{} \\
    \midrule
    \citet{kochProgressMathematicalProgramming2022} & zib03 &19,731,970&29,128,799& 104,422,573& 16.5 (h) & 916\\
    \midrule
      \multirow{4}{*}{Pagerank}  & rand\_1m\_nodes & 1,000,001 & 1,000,000 & 7,999,982 & -    & 3.56  \\
       & rand\_10m\_nodes & 10,000,001 & 10,000,000 & 79,999,982 & -    & 44.22  \\
       % & roadnet-ca & 1,971,282 & 1,971,281 & 9,475,776 & 38.07  & 2.38  \\
       % & cit-patents & 3,774,769 & 3,774,768 & 24,068,483 & -    & 7.71  \\
       & com-livejournal & 3,997,963 & 3,997,962 & 77,358,302 & -    & 21.07  \\
       & soc-livejournal1 & 4,847,572 & 4,847,571 & 78,170,533 & -    & 22.26  \\
       \midrule
      \multirow{2}{*}{Unit Com.}  
      & ds1 & 641,037 & 659,145 & 21,577,566 &  592  &  81\\  
      & ds2 & 641,037 & 659,145 & 21,577,566 &  606  &  108 \\  
       \midrule
      \multirow{4}{*}{Supply-chain}  
      & inv-10 & 4,035,449 & 3,758,458 & 15,264,380   & -    & 1636  \\
      & inv-20 & 8,368,795 & 7,810,584 & 31,718,673   & -    & 1157  \\
      & inv-40 & 13,186,756 & 12,066,105 & 49,528,729 & -    & 6032  \\
      & inv-60 & 16,227,780 & 14,544,689 & 60,372,404 & -    & 11102 \\  
      \midrule
      \multirow{5}{*}{QAP \cite{adams_improved_1994}}  
      & wil50 & \multirow{5}{*}{3,437,600} & \multirow{5}{*}{6,252,500} & \multirow{5}{*}{19,125,000}   & -    & 96  \\
      & lipa50a &  &  &    & -    & 71  \\
      & lipa50b &  &  &  & -    & 66  \\
      & tai50a &  &  &  & -    & 64 \\  
      & tai50b &  &  &  & -    & 437 \\  
      \bottomrule
    \end{tabular}%
    \normalsize
\end{table}%

\section{Conclusions}
In this paper, we introduce an enhanced C implementation of cuPDLP, providing more substantial evidence that GPUs can significantly accelerate first-order linear programming solvers. Specifically, \cupdlp{} demonstrates compelling performance on classic LP benchmarks, including the LP relaxations of MIPLIB 2017 instances and Mittelmann's LP benchmark in finding primal-dual feasible points even compared to commercial IPM-based solvers like COPT. The latter is known to be highly challenging for CPU-based PDLP and other first-order algorithms. Notably, \cupdlp{} outperforms COPT on a handful of instances in these benchmark datasets. According to a comparison to benchmarking results in \cite{mittelmann_lp_2023}, \cupdlp{} is highly competitive as an off-the-shelf linear programming solver.

Furthermore, on selected instances such as PageRank problems, \cupdlp{} exhibits distinct advantages in solving extremely large-scale LP problems, a formidable challenge for traditional IPM and simplex solvers. The results presented in this paper align closely with the conclusions drawn from cuPDLP.jl, emphasizing that an enhanced C implementation, coupled with standard techniques for linear programming, including presolving and scaling, for example, can yield significant improvements. According to the findings of this study, GPUs have a revolutionary impact on traditional LP solvers, and we anticipate that this influence can carry over to other general-purpose solvers for nonlinear programming.

\clearpage
\bibliography{ref}

\begin{thebibliography}{14}
\providecommand{\natexlab}[1]{#1}
\providecommand{\url}[1]{\texttt{#1}}
\expandafter\ifx\csname urlstyle\endcsname\relax
  \providecommand{\doi}[1]{doi: #1}\else
  \providecommand{\doi}{doi: \begingroup \urlstyle{rm}\Url}\fi

\bibitem[Adams and Johnson(1994)]{adams_improved_1994}
Warren~P. Adams and Terri~A. Johnson.
\newblock Improved linear programming based lower bounds for the quadratic
  assignment problem.
\newblock In {Panos M. Pardalos} and {Henry Wolkowicz}, editors,
  \emph{Quadratic {Assignment} and {Related} {Problems}: {Dimacs} {Workshop},
  {May} 20-21, 1993}, number v. 16 in {DIMACS} series in discrete mathematics
  and theoretical computer science. American Mathematical Society, Providence,
  R.I, 1994.
\newblock ISBN 978-0-8218-6607-8.
\newblock Meeting Name: DIMACS Workshop on Quadratic Assignment and Related
  Problems.

\bibitem[Applegate et~al.(2021)Applegate, Hinder, Lu, and
  Lubin]{applegate_faster_2021}
David Applegate, Oliver Hinder, Haihao Lu, and Miles Lubin.
\newblock Faster {First}-{Order} {Primal}-{Dual} {Methods} for {Linear}
  {Programming} using {Restarts} and {Sharpness}.
\newblock \emph{arXiv:2105.12715 [math]}, August 2021.
\newblock URL \url{http://arxiv.org/abs/2105.12715}.
\newblock arXiv: 2105.12715.

\bibitem[Applegate et~al.(2022)Applegate, Díaz, Hinder, Lu, Lubin, O'Donoghue,
  and Schudy]{applegate_practical_2022}
David Applegate, Mateo Díaz, Oliver Hinder, Haihao Lu, Miles Lubin, Brendan
  O'Donoghue, and Warren Schudy.
\newblock Practical {Large}-{Scale} {Linear} {Programming} using
  {Primal}-{Dual} {Hybrid} {Gradient}, January 2022.
\newblock URL \url{http://arxiv.org/abs/2106.04756}.
\newblock arXiv:2106.04756 [math].

\bibitem[Applegate et~al.(2023)Applegate, Hinder, Lu, and
  Lubin]{applegate_faster_2023}
David Applegate, Oliver Hinder, Haihao Lu, and Miles Lubin.
\newblock Faster first-order primal-dual methods for linear programming using
  restarts and sharpness.
\newblock \emph{Mathematical Programming}, 201\penalty0 (1-2):\penalty0
  133--184, September 2023.
\newblock ISSN 0025-5610, 1436-4646.
\newblock \doi{10.1007/s10107-022-01901-9}.
\newblock URL \url{https://link.springer.com/10.1007/s10107-022-01901-9}.

\bibitem[Burkard et~al.(1997)Burkard, Karisch, and Rendl]{burkard1997qaplib}
Rainer~E Burkard, Stefan~E Karisch, and Franz Rendl.
\newblock Qaplib--a quadratic assignment problem library.
\newblock \emph{Journal of Global optimization}, 10:\penalty0 391--403, 1997.

\bibitem[Ge et~al.(2022)Ge, Huangfu, Wang, Wu, and Ye]{ge_cardinal_2022}
Dongdong Ge, Qi~Huangfu, Zizhuo Wang, Jian Wu, and Yinyu Ye.
\newblock Cardinal {Optimizer} ({COPT}) {User} {Guide}, October 2022.
\newblock URL \url{http://arxiv.org/abs/2208.14314}.
\newblock arXiv:2208.14314 [cs, math].

\bibitem[Hinder(2023)]{hinder_worst-case_2023}
Oliver Hinder.
\newblock Worst-case analysis of restarted primal-dual hybrid gradient on
  totally unimodular linear programs, September 2023.
\newblock URL \url{http://arxiv.org/abs/2309.03988}.
\newblock arXiv:2309.03988 [math].

\bibitem[Koch et~al.(2022)Koch, Berthold, Pedersen, and
  Vanaret]{kochProgressMathematicalProgramming2022}
Thorsten Koch, Timo Berthold, Jaap Pedersen, and Charlie Vanaret.
\newblock Progress in mathematical programming solvers from 2001 to 2020.
\newblock \emph{EURO Journal on Computational Optimization}, 10:\penalty0
  100031, 2022.
\newblock ISSN 21924406.
\newblock \doi{10.1016/j.ejco.2022.100031}.
\newblock URL
  \url{https://linkinghub.elsevier.com/retrieve/pii/S2192440622000077}.

\bibitem[Leskovec and Krevl(2014)]{snapnets}
Jure Leskovec and Andrej Krevl.
\newblock {SNAP Datasets}: {Stanford} large network dataset collection.
\newblock \url{http://snap.stanford.edu/data}, June 2014.

\bibitem[Lu and Yang(2023{\natexlab{a}})]{lu_cupdlpjl_2023}
Haihao Lu and Jinwen Yang.
\newblock {cuPDLP}.jl: {A} {GPU} {Implementation} of {Restarted}
  {Primal}-{Dual} {Hybrid} {Gradient} for {Linear} {Programming} in {Julia},
  November 2023{\natexlab{a}}.
\newblock URL \url{http://arxiv.org/abs/2311.12180}.
\newblock arXiv:2311.12180 [math].

\bibitem[Lu and Yang(2023{\natexlab{b}})]{lu_geometry_2023}
Haihao Lu and Jinwen Yang.
\newblock On the {Geometry} and {Refined} {Rate} of {Primal}-{Dual} {Hybrid}
  {Gradient} for {Linear} {Programming}, December 2023{\natexlab{b}}.
\newblock URL \url{http://arxiv.org/abs/2307.03664}.
\newblock arXiv:2307.03664 [math].

\bibitem[Mittelmann(2023)]{mittelmann_lp_2023}
Hans Mittelmann.
\newblock {LP} {Feasibility} {Benchmark} (find a {PD} feasible point), November
  2023.
\newblock URL \url{https://plato.asu.edu/ftp/lpfeas.html}.

\bibitem[Nesterov(2014)]{nesterov2014subgradient}
Yu~Nesterov.
\newblock Subgradient methods for huge-scale optimization problems.
\newblock \emph{Mathematical Programming}, 146\penalty0 (1-2):\penalty0
  275--297, 2014.

\bibitem[Xiong and Freund(2024)]{xiong_computational_2024}
Zikai Xiong and Robert~Michael Freund.
\newblock Computational {Guarantees} for {Restarted} {PDHG} for {LP} based on
  "{Limiting} {Error} {Ratios}" and {LP} {Sharpness}, January 2024.
\newblock URL \url{http://arxiv.org/abs/2312.14774}.
\newblock arXiv:2312.14774 [math].

\end{thebibliography}
\bibliographystyle{plainnat}

\clearpage
\appendix

% \begin{landscape}

\section{Complete Results of Mittelmann's LP Feasible Point Benchmark Set}

\footnotesize
\begin{longtable}{l|r|rr|rr|rr}
  \caption{Complete results of Mittelmann feasible point instances }\label{fig.mittelman.time}                                                                                          \\
  \toprule
  \multirow{2}{*}{Method} & \copt     & \multicolumn{2}{|c}{\cupdlp{} (COPT)} & \multicolumn{2}{|c}{\cupdlp{} (CLP)} & \multicolumn{2}{|c}{\cupdlp}                                     \\
                          & $10^{-8}$ & $10^{-4}$                             & $10^{-8}$                            & $10^{-4}$                    & $10^{-8}$ & $10^{-4}$ & $10^{-8}$ \\
  \midrule
  \endfirsthead

  \toprule
  \multirow{2}{*}{Method} & \copt     & \multicolumn{2}{|c}{\cupdlp{} (COPT)} & \multicolumn{2}{|c}{\cupdlp{} (CLP)} & \multicolumn{2}{|c}{\cupdlp}                                     \\
                          & $10^{-8}$ & $10^{-4}$                             & $10^{-8}$                            & $10^{-4}$                    & $10^{-8}$ & $10^{-4}$ & $10^{-8}$ \\
  \midrule
  \endhead
  \midrule
  \multicolumn{4}{r}{{Continued on next page}}                                                                                                                                          \\
  \midrule
  \endfoot

  \bottomrule
  \endlastfoot
  L1\_sixm250obs          & 0.79      & 6.87                                  & 7.99                                 & 3.27                         & 5.29      & 5.34      & 6.59      \\
  Linf\_520c              & 3.07      & 1253.58                               & 3600.00                              & 11.07                        & 3600.00   & 13.29     & 3600.00   \\
  a2864-99blp             & 0.03      & 0.26                                  & 0.34                                 & 0.35                         & 0.52      & 0.53      & 0.92      \\
  bdry2                   & 6.21      & 3600.00                               & 3600.00                              & 3600.00                      & 3600.00   & 3600.00   & 3600.00   \\
  cont1                   & 0.09      & 238.41                                & 602.34                               & 289.77                       & 1337.23   & 258.78    & 621.90    \\
  cont11                  & 0.19      & 386.36                                & 3419.51                              & 227.04                       & 3600.00   & 497.35    & 2352.79   \\
  datt256                 & 0.79      & 0.17                                  & 0.32                                 & 0.12                         & 0.24      & 0.12      & 0.28      \\
  dlr1                    & 26.58     & 13.36                                 & 3600.00                              & 408.29                       & 3600.00   & 3600.00   & 3600.00   \\
\end{longtable}
\normalsize
% \end{landscape}

% \begin{landscape}
\section{Complete Results of 383 MIPLIB 2017 Instances}
    
\footnotesize
\begin{longtable}{l|r|rr|rr|rr}
\caption{Complete results of 383 MIPLIB collection instances }\label{fig.mip383.time}\\
\toprule
\multirow{2}{*}{Method} & \copt     & \multicolumn{2}{|c}{\cupdlp{} (COPT)} & \multicolumn{2}{|c}{\cupdlp{} (CLP)} & \multicolumn{2}{|c}{\cupdlp}                                     \\
                                & $10^{-8}$ & $10^{-4}$                             & $10^{-8}$                            & $10^{-4}$                    & $10^{-8}$ & $10^{-4}$ & $10^{-8}$ \\
        \midrule
        \endfirsthead

        \toprule
        \multirow{2}{*}{Method} & \copt     & \multicolumn{2}{|c}{\cupdlp{} (COPT)} & \multicolumn{2}{|c}{\cupdlp{} (CLP)} & \multicolumn{2}{|c}{\cupdlp}                                     \\
                                & $10^{-8}$ & $10^{-4}$                             & $10^{-8}$                            & $10^{-4}$                    & $10^{-8}$ & $10^{-4}$ & $10^{-8}$ \\ 
\midrule
\endhead
\midrule
\multicolumn{4}{r}{{Continued on next page}} \\
\midrule
\endfoot

\bottomrule
\endlastfoot
          2club200v15p5scn &    0.04 &    0.69 &    1.41 &    0.23 &    0.36 &    0.24 &    0.39 \\
                    30n20b8 &    0.09 &   19.48 &   19.80 &   19.51 &   19.77 &   26.33 &   27.21 \\
                a2864-99blp &    0.02 &    0.23 &    0.36 &    0.42 &    0.71 &    0.80 &    1.64 \\
       academictimetablebig &    1.98 &    0.24 &    0.45 &    0.44 &    0.78 &    0.42 &    0.80 \\
     academictimetablesmall &    0.30 &    0.26 &    0.36 &    0.20 &    0.30 &    0.20 &    0.30 \\
         adult-max5features &    0.22 &    3.99 &   11.03 &    5.61 &   13.07 &    6.06 &   14.41 \\
          adult-regularized &    0.06 &    0.55 &    0.71 &    0.44 &    0.67 &    1.75 &    2.94 \\
                 allcolor58 &    0.04 &    0.08 &    0.14 &    0.07 &    0.15 &    0.13 &    0.20 \\
\end{longtable}
\normalsize

% \end{landscape}
\end{document}